\documentclass[12pt]{article}

\usepackage{t1enc,amssymb,amsbsy,amsthm,amsmath,amsfonts}

\newcommand{\e}{\rm e}

\newtheorem{theorem}{Theorem}[section]
\newtheorem{lemma}[theorem]{Lemma}
\newtheorem{proposition}[theorem]{Proposition}
\newtheorem{corollary}[theorem]{Corollary}

\newtheorem{remark}[theorem]{Remark}

\numberwithin{equation}{section}

\def\ch{\hbox {\rm cosh\,}}
\def\cth{\hbox {\rm coth\,}}
\def\sh{\hbox {\rm sinh\,}}

\def\ch{\hbox {\rm ch}}
\def\e{\hbox {\rm e}}
\def\sh{\hbox {\rm sh}}
\def\P{{\bf P}}
\def\R{{\bf R}}
\def\R{{\bf R}}
\def\Q{{\bf Q}}
\def\E{{\bf E}}

\def\cE{{\cal E}}
\def\cF{{\cal F}}

\def\cC{{\cal C}}
\def\cG{{\cal G}}

\def\disp{\displaystyle{}}

\def\al{\alpha}

\def\be{\beta}
\def\la{\lambda}

\begin{document}

\author{Paavo Salminen\\{\small Åbo Akademi University}
\\{\small Mathematical Department}
\\{\small FIN-20500 Åbo, Finland} \\{\small email: phsalmin@abo.fi}
\and
Pierre Vallois
\\{\small Universit\'e Henri Poincar\'e}
\\{\small D\'epartement de Math\'ematique}
\\{\small F-54506 Vandoeuvre les Nancy, France}
\\{\small email: vallois@iecn.u-nancy.fr}
}

\title{On maximum increase and decrease of Brownian motion}
\maketitle

\begin{abstract}
The joint distribution of maximum
increase and decrease for Brownian motion up to an independent exponential time
is computed. This is achieved by decomposing the Brownian path at the hitting times
of the infimum and the supremum before the exponential time.
It is seen that an important element in our formula is the distribution of
the maximum decrease for the three dimensional Bessel process with drift started
from 0 and stopped at the first hitting of a given level. From the joint distribution
of the maximum
increase and decrease it is possible to calculate the correlation coefficient  between these  at a fixed time and this is seen to be -0.47936... .
\\
\\
\\

Dans cet article nous d\'eterminons la loi conjointe de la plus
grande mont\'ee et de la plus grande descente   d'un  mouvement
brownien arr\^et\'e en un temps exponentiel independant. La preuve
repose sur la d\'ecomposition de la trajectoire brownienne aux
instants o\`u le processus atteint son maximum, resp.   son
minimum, avant le temps exponentiel. La loi de la plus grande
descente d'un processus de Bessel, de dimension trois, issu de $0$
et arr\^et\'e lorsqu'il atteint un niveau fix\'e, joue \'egalement
un r\^ole important. Le coefficient de corr\'elation lin\'eaire de
la grande mont\'ee et de la plus grande descente   d'un  mouvement
brownien arr\^et\'e en temps fixe est d\'etermin\'e :
-0.47936... .
\end{abstract}

\eject
\section{Introduction and notation}
\label{intro}

{\bf 1.} In this paper we are interested in the joint distribution of
the maximum increase and
decrease for a standard Brownian motion, for short BM. 
Let us start with some notation. Let $\Omega:=\cC(\R_+,\R)$ be the space of continuous functions
$\omega:\R_+\mapsto \R$ 
and $X_t(\omega)=\omega(t),\, t\geq 0,$ the coordinate
mappings. With every $\omega$ we associate its lifetime $\zeta(\omega)\in
(0,\infty]$ and consider $X_t$ to be defined for $t<\zeta(\omega).$ 
The standard notation $\cF_t$ is used for the
$\sigma$-algebra generated by the coordinate mappings up
to time $t,$ and we set $\cF:=\cF_\infty.$

Further, $\P^\mu_x $ and $\E^\mu_x $ denote the probability measure and
the expectation operator on $(\Omega,\cF)$ under which the coordinate
process $X=\{X_t\,:\, t\geq 0\}$ is a Brownian motion with drift
$\mu$ started from $x,$  for short BM($\mu$). For simplicity, $\P_x $ and $\E_x $ stand
for the corresponding objects for BM.

The maximum increase before time $t$ is defined as
\begin{equation}
\label{new+}
D^+_t:=\sup_{0\leq u\leq v\leq t}\left(X_v -X_u\right)
\end{equation}
and, analogously, the maximum decrease
\begin{equation}
\label{new-}
D^-_t:=\sup_{0\leq u\leq v\leq t}\left(X_u -X_v\right)
\end{equation}
Notice that, e.g.,
\begin{equation}
\label{e00}
D^-_t=\sup_{0\leq v\leq t}\left(\sup_{0\leq u\leq v}X_u
-X_v\right).
\end{equation}
Using the L\'evy isomorphism, i.e., under $\P_0$
\begin{equation}
\label{e002}
\{X^+_v:=\sup_{0\leq u\leq v}X_u
-X_v\,:\, v\geq 0\}
\ {\stackrel {(d)}{=}}\ 
\{|X_v|\,:\, v\geq 0\},
\end{equation}
where ${\stackrel {(d)}{=}}$ means ``to be identical in law with'',
it follows from identity (\ref{e00}) that
\begin{equation}
\label{en001}
\P_0(D^-_t>a)=\P_0(H_a\wedge H_{-a}<t)
\end{equation}
with 
$$
H_b:=\inf\{t\,:\,X_t=b\},\quad b\in\R,
$$
the first hitting time of $b$ (in the canonical setting)
with the usual convention that $H_b=+\infty$ if the set in the braces is empty.
From the equality (\ref{en001}) applying, e.g.,
\cite{borodinsalminen02} 1.3.0.2 p. 212, and  3.1.1.4 p. 333 
we obtain
\begin{eqnarray*}
&&\hskip-.8cm
\P_0(D^-_t>a)=\sum_ {k=-\infty}^{+\infty} (-1)^k
\int_0^tds\,\frac{(2k+1)a}{\sqrt{2\pi}s^{3/2}}\, {\rm e}^{-(2k+1)^2a^2/2s}
\\
&&\hskip1.5cm
=1- \frac{1}{\sqrt{2\pi t}} \sum_ {k=-\infty}^{+\infty} (-1)^k
\int_{-a}^{a} dx\,\left({\rm e}^{-(x+4ka)^2/2t}  -{\rm e}^{-(x+4ka+2a)^2/2t}\right).
\end{eqnarray*}
Notice also that due to the symmetry of standard Brownian motion we have
$$
\{D^-_t\,:\, t\geq 0\}\ {\mathop=^{\rm{(d)}}}\ \{D^+_t\,:\, t\geq 0\},
$$
and this holds under $\P_x$ for any $x\in\R.$ We refer to Douady, et al. \cite{douady99}
for results concerning the distribution of maximum increase and related functionals up to a fixed time in the case of Brownian motion.

For Brownian motion with drift the result corresponding  to
(\ref{e002}) states that under $\P^\mu_x $ 
the process $\{X^+_v\,:\, v\geq 0\}$ 
is a reflected Brownian motion on $\R_+$ with drift $-\mu,$ for short RBM($-\mu$), (see, e.g., Harrison \cite{harrison85} p. 49, and  McKean
\cite{mckean69} p. 71), more precisely, it is a
diffusion on $\R_+$ with basic characteristics as given in
\cite{borodinsalminen02} A1.16 p. 129.  The
probability measure on $(\Omega,\cF)$ associated with $X ^+$ (under  $\P^{\mu}$)
is denoted by $\P^{-\mu,+}.$  Clearly, we have
for a given $a>0$
\begin{equation}
\label{e1}
\P^{\mu}_0(D^-_t>a)=\P^{-\mu,+}_0(H_a<t).
\end{equation}
Similarly, defining 
$$
X^-_v:=X_v-\inf_{0\leq u\leq v}X_u
$$
it holds under $\P^\mu$ that the process $X^-$ is a RBM($\mu$). Letting $\P^{\mu,+}$ denote the measure associated
with $X^-$ we have 
\begin{equation}
\label{ee1}
\P^{\mu}_0(D^+_t>a)=\P^{\mu,+}_0(H_a<t).
\end{equation}
For an explicit expression of the  $\P^{\mu}$-distribution of $D^-_t$, see Domin\'e
\cite{domine96} where the method based on spectral representations is
used. In Magdon-Ismail  et.al.  \cite{magdonismail04}
formulas for the mean of $D^-_t$ are derived. 

\medskip
\noindent
{\bf 2.} Unfortunately we are not able to determine explicitly the distribution
of  $(D^+_t,D^-_t),$ but replacing $t$ by $T,$ that is, an exponentially distributed
random variable independent of $X$ with mean $1/\lambda,$ 
allows us to find the $\P$-distribution of $(D^+_T,D^-_T),$
see Proposition \ref{cor43} and \ref{cor44}. We remark that the marginal $\P^\mu$-distributions
of  $D^+_T$ and $D^-_T$ are easily computed from (\ref{ee1}) and
(\ref{e1}), respectively. Indeed, using standard diffusion theory and some explicit formulas (see e.g.
\cite{borodinsalminen02} p. 18 and 129) yield
\begin{equation}
\label{en01}
\P^{\mu}_0(D^-_T>a)=\E^{-\mu,+}_0\left(\exp(-\lambda
H_a)\right)=1/\psi_\lambda(a;-\mu)
\end{equation}
and 
\begin{equation}
\label{en02}
\P^{\mu}_0(D^+_T>a)=1/\psi_\lambda(a;\mu)
\end{equation}
with
$$
\psi_\lambda(a;\nu):=
{\rm e}^{-\nu a}\left(\ch(a\,\sqrt{2\lambda+\nu^2})+\frac{\nu}{\sqrt{2\lambda+\nu^2}}
\sh(a\,\sqrt{2\lambda+\nu^2})\right).
$$

In our approach for finding the joint distribution we consider first 
the case where the infimum is attained before the supremum. 
In this case it is clear that the maximum increase is nothing but the
difference of the supremum and the infimum, and, in a sense, we have reduced
the problem to the problem for finding the distribution of the maximum
decrease. 
The opposite
case  where the supremum is attained before the infimum is clearly treated using symmetry. 

It is natural when the infimum is attained before the supremum to decompose the exponentially stopped Brownian path into three
parts: 
\begin{description}
\item{\hskip.3cm$\circ$\ } the first part is up to the hitting time of the infimum,
\item{\hskip.3cm$\circ$\ } the second part is from the hitting time of the infimum to the hitting
time of the supremum 
\item{\hskip.3cm$\circ$\ } the third part is from the hitting time of the supremum to the exponential time.
\end{description}

\noindent
We prove in Theorem \ref{deco2} that these three parts are
conditionally independent given the infimum and the supremum, and
find their distributions in terms of the three dimensional Bessel
processes with drift. Our approah is mainly based on the  $h$-transform
techniques, excursion theory and path decompositions of
Brownian motion with drift.

The above described path decomposition up to $T$ permits us to
determine the joint distribution of $(D^+_T,D^-_T)$ since now $D^+_T=S_T-I_T$ and, under this
decomposition, $D^-_T$ is the maximum of the maximum decreases of the
three conditionally independent fragments.    
To find the distribution of the maximum decrease for the
first and the third part is fairly straightforward diffusion theory. 
To compute the maximum decrease for the second part 
is equivalent for finding the distribution of the maximum decrease for a three dimensional Bessel
process with drift (see Proposition \ref{prop2}). 


Although the distribution and the density function of
$(D^+_T,D^-_T)$ are complicated it is possible to determine by the scaling property of BM
the covariance between $D^+_t$ and $D^-_t$ and this is given by 
$$
\E(D^+_t\,D^-_t)= (1-2\log 2 +2\beta(2))\, t,
$$
where 
$$
\beta(2):= \sum_{k=0}^\infty(-1)^k\,(2k+1)^{-2}=0.91596...
$$
is Catalan's constant. Hence, the correlation coefficient $\rho$ 
between  $D^+_T$ and $D^-_T$ is easily obtained to  be
$$
\rho:=\frac{\E(D^+_t\,D^-_t)-(\E(D^+_t))^2}{{\bf Var}(D^+_t)}
=-0.47936... .
$$

\noindent
{\bf 3.} One motivation to study
the maximum decrease and increase comes from mathematical finance where
the maximum decrease, also called maximum drawdown (MDD), is used to
quantify the riskyness of a stock or any other asset. Related
measures used hereby are e.g. the recovery time from MDD and the duration of MDD.
Our interest to the problem discussed in the paper
was arised by Gabor Szekely who asked for an expression for the covariance
between $D^+_t$ and $D^-_t.$

\noindent
{\bf 4.} The paper is organised so that in the next section we find the
distribution of the maximum decrease of a stopped Brownian motion with
positive drift. In fact, we compute this distribution under the restriction
that the process does not hit some negative level, and proceed from
here to
the distribution of the maximum decrease for a three dimensional Bessel
process with drift. In the third section path decompositions are
discussed. To prove our main path decomposition Theorem \ref{deco2}, we first prove a
decomposition of the Brownian trajectory $\{B_t\,:\, t\leq T\}$
conditionally on $I_T$ (see Theorem \ref{deco1}). The fourth section
is devoted to computation and analysis of the law of  $(D^-_T,D^+_T).$

\section{Maximum decrease for stopped Brownian motion with drift}
\label{sec1}

According to (\ref{new-}), the maximum decrease up to the first hitting time
of a given level $\beta$ is
$$
D^-_{H_\beta}:=
\sup\{X_u
-X_v\,:\, 0\leq u\leq v\leq H_\beta\}.
$$
In this section we consider the $\P^\mu$-distribution of $D^-_{H_\beta}$ 
under some additional conditions and conditioning.
Recall that  
\begin{equation}
\label{e110}
S^{\mu}(x):=\frac{1}{2\mu}(1-{\rm e}^{-2\mu x})
\end{equation}
is the scale function of BM($\mu)$
and for $a<x<b$
\begin{equation}
\label{e112}
\P^{\mu}_x(H_{a}<H_{b})=\disp{\frac{S^{\mu}(b)-S^{\mu}(x)}{S^{\mu}(b)-S^{\mu}(a)}}.
\end{equation}

\begin{proposition}
\label{prop1}
{\sl For all nonnegative $\al, \be,$ and $u$ 
$$
\P^{\mu}_0(D^-_{H_\beta}<u\,,\,H_\be<H_{-\al})
=
\begin{cases}
\exp\left(-\disp{\frac{\be}{S^{-\mu}(u)}}\right),&\hskip-2cm u\leq \al,\\
&\\
\disp{\frac{S^{\mu}(\al)}{S^{\mu}(u)}}
\exp\left(-\disp{\frac{\be+\al-u}{S^{-\mu}(u)}}\right), & \\
&\hskip-2cm \al\leq u\leq \al+\be,\\
\disp{\frac{S^{\mu}(\al)}{S^{\mu}(\al+\be)}},&\hskip-2cm \al+\be\leq u,\\
\end{cases}
$$ 
In particular, 
$$
\P^{\mu}_0(D^-_{H_\beta}<u)
=
\exp\left(-\disp{\frac{\be}{S^{\mu}(u)}}\right).
$$
For standard Brownian motion, i.e., $\mu=0,$ the above formulas hold
with $S^{0}(u)=u.$ 
}
\end{proposition}
\begin{proof} 
We assume that 
$X_0=0,$ and define for $a>0$
$$
H_{a+}:=\inf\{t\,:\,X_t>a\}.
$$ 
For a given $a>0,$ in the case $H_{a+}>H_{a},$ let
$$
\Xi^+_a(u):=a-X_{u+H_{a}},\quad  0\leq u<H_{a+}
  -H_{a},
$$
and, if  $H_{a+}= H_{a},$ take $\Xi^+_a:=\partial$,
where $\partial$ is some fictious (cemetary) state.
The process 
$$
\Xi^+=\{ \Xi^+_a\,:\, a\geq 0\}
$$
is called the excursion process, associated with $X,$ for excursions
under the running maximum. Let, further,
$$
M_a:=\sup\{\Xi^+_a(u)\,:\, 0\leq u<H_{a+} - H_{a}\}.
$$
An obvious but important fact is that
\begin{equation}
\label{e11}
D^-_{H_\beta}=\sup_{a<\beta}M_a.
\end{equation}
Introduce also for $u>0$
$$
\xi_u:=\inf\{a\geq 0\,:\, M_a>u\},
$$
and 
$$
\xi^\circ_u:=\inf\{a\geq 0\,:\, a-M_a<-u\}.
$$
Then it holds for positive $\al,\be,$ and $u$
$$
\{D^-_{H_\beta}\leq u\}=\{\forall\ a\in(0,\beta)\,:\, M_a\leq u\}=\{\xi_u\geq\be\},
$$
and 
$$
\{H_\be<H_{-\al}\}=\{\forall\ a\in(0,\beta)\,:\,a- M_a>-\al\}=\{\xi^\circ_\al>\be\};
$$
hence, for $0<u<\al+\be$
\begin{equation}
\label{e121}
\P^{\mu}_0(D^-_{H_\beta}<u\,,\,H_\be<H_{-\al})=\P^{\mu}_0(\xi_u> \be\,,\,\xi^\circ_\al> \be).
\end{equation}
Since $X^+$ under $\P^{\mu}$ is identical in law with RBM($-\mu$) it
follows that the excursion process $\Xi^+$ is
identical in law with the usual excursion process of RBM($-\mu$) for
excursions from 0 to 0. Consequently, see Pitman and
Yor \cite{pitmanyor82}, 
$$
\Pi=\{(a,M_a)\,:\,a\geq 0\}
$$ 
is a homogeneous Poisson point process with
the characteristic measure
$$
\nu(da,dm)=da\,n(dm),
$$
where for $m>0$
$$ 
n((m,+\infty))=1/ S^{-\mu}(m).
$$
Introduce the sets  
$$
A:=[0,\be)\times[u,+\infty)\quad {\rm and}\quad B:=\{(a,m)\,:\, 0\leq
    a<\be\,,\, a-m<-\al\},
$$
and let $N$ denote the counting measure associated with $\Pi.$ Now we
have 
\begin{eqnarray*}
&&
\P^{\mu}_0(\xi_u> \be\,,\,\xi^\circ_\al> \be)=\P^{\mu}_0(N(A\cup B)=0)\\
&&
\hskip3.7cm
=\exp\left(-\nu(A\cup B)\right).
\end{eqnarray*}
It is straightforward to compute $\nu(A\cup B)$ for different values
on $\al,\be,$ and $u,$ and we leave this to the reader. Consequently,
by (\ref{e121}), the claimed formula is obtained.
\end{proof}

Obviously, $D^+_{H_\be}\geq \be$ and for $z>0$
$$
\{D^+_{H_\be}-\be<z\}=\{H_\be<H_{-z}\}.
$$
Consequently, we have from Proposition \ref{prop1} the following
corollary giving an expression for the joint
distribution of  $D^+_{H_\be}$ and $D^-_{H_\be}.$ Notice that 
$$
 D^-_{H_\be} \leq D^+_{H_\be}\leq  D^-_{H_\be} +\be
$$
explaining the three cases below.

\begin{corollary}
\label{cor1} For $v\geq \be$ and $\mu\geq 0$ 
$$
\P^{\mu}_0(D^-_{H_\beta}<u\,,\,D^+_{H_\beta}<v)
=
\begin{cases}
\exp\left(-\disp{\frac{\be}{S^{-\mu}(u)}}\right),&\hskip-2cm u\leq v-\be,\\
&\\
\disp{\frac{S^{\mu}(v-\be)}{S^{\mu}(u)}}
\exp\left(-\disp{\frac{v-u}{S^{-\mu}(u)}}\right), & \\
&\hskip-2cm v-\be\leq u\leq v,\\
\disp{\frac{S^{\mu}(v-\be)}{S^{\mu}(v)}},&\hskip-2cm v\leq u,\\
\end{cases}
$$ 
For standard Brownian motion, i.e., $\mu=0,$ the above formulas hold
with $S^{0}(u)=u.$ 
\end{corollary}

We proceed by developing the result in Proposition \ref{prop1} for a 
3-dimensional Bessel process with drift $\mu>0$, for short
BES($3,\mu$). We recall that BES($3,\mu$) is a linear
diffusion with the generator
\begin{equation}
\label{e13}
\cG^{R,\mu}=\frac 12\, \frac{d^2}{dx^2}+\mu\coth(\mu x)\,\frac{d}{dx},\quad x>0.
\end{equation}
The notation $\Q^\mu_x$ is used for the probability measure on 
the canonical space $\Omega$ associated with BES($3,\mu$) when started
from $x\geq 0.$ In the case $\mu=0$ the corresponding measure is simply denoted by  
$\Q_x$ and the generator is given by 
\begin{equation}
\label{e131}
\cG^{R}=\frac 12\, \frac{d^2}{dx^2}+\frac 1x\,\frac{d}{dx},\quad x>0.
\end{equation}

The following lemma is a fairly well known example on
$h-$transforms. To make the presentation more self contained we give a
short proof. It is also interesting to compare the result with Lemma
\ref{lem2} in the next section.
\begin{lemma}
\label{lem1}
{\sl
Let $0<x<y$ be given. The Brownian motion with drift $\mu$ started from $x>0,$ killed
at the first hitting time of $y,$
and conditioned to hit $y$ before 0 is identical in law with
a 3-dimensional Bessel process with drift $|\mu|$ started from $x$ and
killed at the first hitting time of $y.$
} 
\end{lemma}
\begin{proof} 
In our canonical space of continuous functions with $X_0=x$ we have 
for a given $t>0$
$$
\{t<H_y<H_0\}= \{t<H_y\wedge H_0\,,\,
H_y\circ\theta_t<H_0\circ\theta_t \},
$$
where $\theta_\cdot$ is the usual shift operator, i.e., $X_s\circ
\theta_t=X_{s+t}.$ Hence, for any 
$A_t\in \cF_t,$  
\begin{eqnarray*}
&&\P^{\mu}_x(A_t\, ,\, t<H_y\,| \,H_y<H_{0})
\\
&&\hskip2cm
=
\P^{\mu}_x(A_t\, ,\, t<H_y\wedge H_0\,,
\,H_y\circ\theta_t<H_0\circ\theta_t)/\P^{\mu}_x(H_y<H_{0})
\\
&&\hskip2cm
=\E^{\mu}_x(h_1(X_t)\,;\, A_t\, ,\, t<H_y\wedge H_0)/h_1(x)
\end{eqnarray*}
by the Markov property, where
$$
h_1(x):=\P^{\mu}_x(H_y<H_{0} )= 
\disp{\frac{S^{\mu}(x)-S^{\mu}(0)}{S^{\mu}(y)-S^{\mu}(0)}}
=
\disp{\frac{1-{\rm e}^{-2\mu x}}{1-{\rm
      e}^{-2\mu y}}}.
$$
Consequently, the desired conditioning can be realized by taking 
the Doob $h-$transform (with $h=h_1$) of the Brownian motion with
drift $\mu$ killed at
$H_{y}\wedge H_{0}.$ 
The generator of the $h-$transform is 
$$
\cG^{h_1}:=
\frac 12\, \frac{d^2}{dx^2}+\mu\,\frac{d}{dx} +\frac{h'(x)}{h(x)}\,\frac{d}{dx},
$$
which is easily seen to coincide with (\ref{e13}) with $|\mu|$ instead
of $\mu.$
\end{proof}

\begin{remark}
\label{rem22}
{\sl Analogously as above, it can be proved that BM($\mu$) started from $x>0$ and conditioned not to hit 0
is identical in law with BES(3,$|\mu|$) started from $x.$ 
}
\end{remark}

\begin{proposition}
\label{prop2}
{\sl
For $\be >u>0$ 
$$
\Q^{\mu}_0(D^-_{H_\be}<u)=\disp{\frac{S^{-\mu}(\be)}{S^{-\mu}(u)}}
\exp\left(-\disp{\frac{\be-u}{S^{-\mu}(u)}}-2\mu(\be-u)\right).
$$
For the 3-dimensional Bessel process without drift, i.e., $\mu=0,$ the
above formula holds with $S^{0}(u)=u.$
}
\end{proposition}
\begin{proof}
Note that under $\Q^{\mu}_0,$ it holds a.s. on $\{D^-_{H_\be}>u\}$ that
$$
D^-_{H_\be}=D^-_{H_\be}\circ \theta_{H_u}.
$$
Therefore, applying the strong Markov property at time $H_u$ yields
$$
\Q^{\mu}_0(D^-_{H_\be}>u)=\Q^{\mu}_u(D^-_{H_\be}>u),
$$
and from Lemma \ref{lem1}
\begin{eqnarray*}
&&
\Q^{\mu}_u(D^-_{H_\be}>u)=\P^{\mu}_u(D^-_{H_\be}>u\, |\,H_\be<H_0 )
\\
&&\hskip2.75cm
=\P^{\mu}_u(D^-_{H_\be}>u\, ,\,H_\be<H_0 )/\P^{\mu}_u(H_\be<H_0 )
\\
&&\hskip2.75cm
=\P^{\mu}_0(D^-_{H_{\be-u}}>u\, ,\,H_{\be-u}<H_{-u} )/\P^{\mu}_0(H_{\be-u}<H_{-u} )
.
\end{eqnarray*}
The proof is
now easily completed from Proposition \ref{prop1}. 
\end{proof}

\section{Path decompositions}
\label{sec2}

The main path decomposition results presented in this section are stated in Theorems \ref{deco1} 
 and \ref{deco2}. In the first one we consider the decomposition at the global infimum and the second one gives, roughly speaking, the decomposition of the post part of the  previous decomposition at its global supremum.

We begin with by stating the following lemma which is proved similarly as Lemma \ref{lem1}.
\begin{lemma}
\label{lem2}
{\sl
Let $0<x<y$ be given. The Brownian motion started from $x>0,$ killed
at the first hitting time of $y,$
and conditioned by the event $H_y<H_0\wedge T,$ 
where $T$ is an exponentially  distributed
random variable with parameter $\lambda$ independent of the Brownian motion,
is identical in law with
BES(3,$\sqrt{2\lambda}$) started from $x$ and
killed at the first hitting time of $y.$
} 
\end{lemma}
\begin{proof} 
We adapt the proof of Lemma \ref{lem1} to our new situation. For $t>0$
we have 
$$
\{t<H_y<H_0\wedge T\}= \{t<H_y\wedge H_0\wedge T\,,\,
H_y\circ\theta_t<(H_0\wedge T)\circ \theta_t \}.
$$
Hence, using the memoryless property of $T$ we get for $A_t\in\cF_t$
\begin{eqnarray*}
&&\P_x(A_t\, ,\, t<H_y\,| \,H_y<H_{0}\wedge T )
\\
&&\hskip2cm
=\E_x(h_2(X_t)\,;\, A_t\, ,\,
t<H_y\wedge H_0\wedge T)/h_2(x)
\end{eqnarray*}
with
\begin{eqnarray*}
&&
h_2(x)=\P_x(H_y<H_{0}\wedge T)
= \P_x(H_y<H_{0}\,,\,H_y< T) 
\\
&&
\hskip1.1cm
={\rm sh}(x\sqrt{2\lambda})/{\rm sh}(y\sqrt{2\lambda})
\end{eqnarray*}
(see \cite{borodinsalminen02} 1.3.0.5(b) p.212). 

Consequently, the desired conditioning can be realized by taking 
the Doob $h-$transform (with $h=h_2$) of a Brownian motion killed at time 
$H_{y}\wedge H_{0}\wedge T.$ The generator of the $h-$transform can 
be computed in the usual way, and is seen to coincide with the
generator of BES($3,\sqrt{2\lambda}$) (see (\ref{e13})).
\end{proof} 

We let, throughout the paper, $T$ denote an exponentially with parameter $\lambda$
distributed random variable independent of $X$ under $\P_0,$ and define 
$$
I_T:=\inf\{X_t:\, 0\leq t\leq T\}\quad,\quad S_T:=\sup\{X_t:\, 0\leq t\leq T\}
$$
and
$$
H_I:=\inf\{t:\,  X_t = I_T\}\quad ,\quad H_S:=\inf\{t:\,  X_t = S_T\}.
$$
Next we discuss the path decomposition at the global infimum
for Brownian motion killed at $T.$ If nothing else is stated the
coordinate process is considered under $\P_0.$

\begin{theorem}
\label{deco1}
{\bf 1.} 
The processes 
$\{X_{t}:\, 0\leq t< H_I\}$ and  $\{X_{T-t}-X_T:\, 0\leq t< T-H_I\}$
are independent and identically distributed.

\noindent
{\bf 2.} Given $I_T=a$ 
\begin{description}
\item{1.} the pre-$H_I$ process $\{X_t:\, 0\leq t< H_I\}$
and the post-$H_I$ process $\{X_{H_I+t}:\, 0\leq t< T-H_I\}$
are independent.
\item{2.} the pre-$H_I$ process is identical
in law with a BM($-\sqrt{2\lambda}$) killed when it hits $a,$
\item {3.}  the post-$H_I$ process
is identical in law with the diffusion $Z$ started from $a$ and having the generator
\begin{equation}
\label{h0}
{\cal G}^Zu(x)=\frac 12 u''(x) +\frac{h'_3(x-a)}{h_3(x-a)}u'(x)-\frac{\lambda}{h_3(x-a)}u(x),
\end{equation}
where $x>a$ and
\begin{equation}
\label{h1}
h_3(y):=\P_y(T<H_0)=1-{\rm e}^{-y\,\sqrt{2\lambda}},\quad y>0.
\end{equation}
Moreover, $Z$ is the Doob $h$-transform with $h=h_3(\cdot-a)$ of BM killed at time $T\wedge H_a.$
\end{description}
\end{theorem}

\begin{proof} {\sl a)} We have two different proofs: the first one is "direct" in a sense that we compute the conditional finite dimensional distributions, the second one relies on excursion theory of Brownian motion. From our point of view the both proofs contain interesting elements which motivates the presentation of both of these. 
\hfill\break\hfill
{\sl b)} We begin with the {\sl direct proof} of claim 2. Define for $s<t$
$$
I_{s,t}:=\inf\{X_u:\, s\leq u\leq t\},\quad  I_t:=I_{0,t}
$$
and
$$
H_{I_{s,t}}:=\inf\{u\in(s,t):\,  X_u = I_{s,t}\}, \quad  H_{I_t}:= H_{I_{0,t}}.
$$
Let $u,v,$ and $t$ be given such that $0<u<v<t.$
For positive integers $n$ and $m$ introduce
$
0<u_1<\dots<u_n<u
$
and 
$ 0<v_1<\dots<v_m
$
with $v_m+v<t.$
Define also 
\begin{equation}
\label{An}
A_n:=\{X_{u_1}\in dx_1,\dots,X_{u_n}\in dx_n\},
\end{equation}
and
\begin{equation}
\label{Bm}
B_{m}:=\{X_{v_1}\in dy_1,\dots,X_{v_m}\in dy_m\}.
\end{equation}
Consider now for $u<s<v$ 
\begin{eqnarray*}
&&
\hskip-.5cm
\P_0(A_n, I_t\in da, H_{I_t}\in ds, B_{m}\circ\theta_v, X_t\in dz) \\
&&
\hskip-.5cm
=
\P_0(A_n, I_{u}>a, I_{u,v}\in da, H_{I_{u,v}}\in ds,
I_{v,t}>a, B_{m}\circ\theta_v, X_t\in dz)\\
&&
\hskip-.5cm
=
\P_0\left(A_n, I_{u}>a; \P_0\left( I_{u,v}\in da, H_{I_{u,v}}\in ds,
I_{v,t}>a, B_{m}\circ\theta_v, X_t\in dz\,|\, \cF_u\right)\right).
\end{eqnarray*}
Further,
\begin{eqnarray*}
&&
\hskip-.5cm
\P_0\left( I_{u,v}\in da, H_{I_{u,v}}\in ds,
I_{v,t}>a, B_{m}\circ\theta_v, X_t\in dz\,|\, \cF_u\right)\\
&&
\hskip.5cm
=
\P_0\left( I_{u,v}\in da, H_{I_{u,v}}\in ds;
\P_0\left(I_{v,t}>a, B_{m}\circ\theta_v, X_t\in dz\,|\, \cF_v\right)\,|\, \cF_u\right)\\
&&
\hskip.5cm
=
\P_0\left( I_{u,v}\in da, H_{I_{u,v}}\in ds;
\P_0\left(I_{v,t}>a, B_{m}\circ\theta_v, X_t\in dz\,|\, X_v\right)\,|\, X_u\right)
\end{eqnarray*}
by the Markov property. Letting $p^+$ denote the transition density (with respect to
$2\,dx$) of BM killed when
it hits $a$ and
writing $X_v=z_2$ we have 
\begin{eqnarray*}
&&\hskip-1.5cm
\P_0\left(I_{v,t}>a, B_{m}\circ\theta_v, X_t\in dz\,|\, X_v\right)\\
&&\hskip1cm
=p^+(v_1;z_2,y_1)\, 2dy_1\cdot\dots\cdot
p^+(t-v-v_m;y_{m},z)\,2dz\\ 
&&\hskip1cm
=: F_{v_1,\dots,v_m,t-v}(z_2,y_1,\dots,y_m, z)\,2dy_1\dots 2dz.
\end{eqnarray*}
Introduce 
\begin{equation}
\label{b01}
\eta_x(a,\al):=\P_x(H_a\in d\al)/d\al,
\end{equation}
and recall the formula due to L\'evy  
\begin{eqnarray}
\label{levy}
&&
\hskip-1cm
\nonumber
\P_x(I_\be\in da, H_{I_\be}\in d\al, X_\be\in dz)\\
&&\hskip2cm
=\eta_x(a,\al)\,
\eta_z(a,\be-\al)\,d\al\, 2dz\,da,\quad \al<\be.
\end{eqnarray}
Applying (\ref{levy}) and putting $X_u=z_1$ we obtain 
\begin{eqnarray*}
&&\hskip-1cm
\P_0\left( I_{u,v}\in da, H_{I_{u,v}}\in ds;
\P_0\left(I_{v,t}>a, B_{m}\circ\theta_v, X_t\in dz\,|\, \cF_v\right)\,|\, \cF_u\right)\\
&&
=
\int_a^\infty 2dz_2\, \eta_{z_1}(a,s-u)\,
\eta_{z_2}(a,v-s)\,da\,ds\\
&&
\hskip2cm
\times 
 F_{v_1,\dots,v_m,t-v}(z_2,y_1,\dots,y_m, z)\,2dy_1\dots 2dy_m\,2dz,
\end{eqnarray*}
and, finally,
\begin{eqnarray}
\label{b1}
&&
\hskip-.5cm
\P_0(A_n, I_t\in da, H_{I_t}\in ds, B_{m}\circ\theta_v, X_t\in dz) \\
&&
\nonumber
\hskip-.5cm
=\int_a^\infty 2\,dz_1\int_a^\infty 2\,dz_2\
p^+(u_1;0,x_1)\, 2dx_1\cdot\dots\cdot
p^+(u-u_n;x_{n},z_1)\\
&&
\nonumber
\hskip1cm
\times 
 \eta_{z_1}(a,s-u)\,
\eta_{z_2}(a,v-s)\,da\,ds
\\
&&
\nonumber
\hskip1.5cm
\times 
 F_{v_1,\dots,v_m,t-v}(z_2,y_1,\dots,y_m, z)\,2dy_1\dots 2dy_m\,2dz\\
&&
\nonumber
\hskip-.5cm
=\int_a^\infty 2\,dz_1\int_a^\infty 2\,dz_2\, 
 F_{u_1,\dots,u_n,u}(0,x_1,\dots,x_n, z_1)\,2dx_1\dots 2dx_n
\\
&&
\nonumber
\hskip1cm
\times 
 \eta_{z_1}(a,s-u)\,
\eta_{z_2}(a,v-s)\,da\,ds\,
\\
&&
\nonumber
\hskip1.5cm
\times 
 F_{v_1,\dots,v_m,t-v}(z_2,y_1,\dots,y_m, z)\,2dy_1\dots 2dy_m\,2dz.
\end{eqnarray}
Replacing in (\ref{levy}) the deterministic time $\be$ with the
exponential time $T$ yields for  $a<0$ and $a<z$
\begin{eqnarray}
\label{levy2}
&&
\nonumber
\P_0(I_T\in da, H_{I}\in ds, X_T\in dz)\\
&&
\hskip3cm
=\eta_0(a,s)\,\lambda\, {\rm
  e}^{-\lambda s}\ \E_z\left({\rm
  e}^{-\lambda H_a}\right)\,ds\, 2dz\,da,
\end{eqnarray}
and, further,
\begin{equation}
\label{levy3}
\P_0(I_T\in da, H_{I}\in ds)=\sqrt{2\lambda}\, {\rm
  e}^{-\lambda s} \eta_0(a,s)\,da\,ds.
\end{equation}
We operate similarly in (\ref{b1}), i.e., introduce the exponential
time $T$ in place of $t.$ After this we 
integrate over $z,$ and
divide with the expression on the r.h.s. in (\ref{levy3}) and obtain 
for $u<s<v$
\begin{eqnarray}
\label{b2}
&&\hskip-1cm
\nonumber
\P_0(A_n,  B_{m}\circ\theta_v\,|\,  I_T=a, H_{I}=s)\\
&&\hskip1cm
= \widehat F_{u_1,\dots,u_n}(x_1,\dots,x_n; a,s)\,2dx_1\dots 2dx_n\\
&&\hskip2cm\nonumber
\times\ \widehat G_{v_1,\dots,v_m}(y_1,\dots,y_m; a;s,v)\,2dy_1\dots 2dy_m,
\end{eqnarray}
with
\begin{eqnarray}
\label{F}
&&\hskip-1cm
\widehat F_{u_1,\dots,u_n}(x_1,\dots,x_n; a;s)\\
&&
\nonumber
=
\widehat p^+(u_1;0,x_1)\, \cdot\dots\cdot
\widehat p^+(u_n-u_{n-1};x_{n-1},x_n)\
\frac {\widehat \eta_{x_n}(a,s-u_n)}{\widehat \eta_0(a,s)}
\end{eqnarray}
and
\begin{eqnarray}
\label{G}
&&
\widehat G_{v_1,\dots,v_n}(y_1,\dots,y_n; a;s,v)
=\frac {\widehat \eta_{y_1}(a,v-s+v_1)}{\sqrt{2\lambda}}
\,
\\
\nonumber
&&\hskip2.5cm
\times \widehat p^+(v_2-v_1;y_1,y_2)\cdot\dots\cdot
\widehat p^+(v_m-v_{m-1};y_{m-1},y_m) \,h_3(y_m-a),
\end{eqnarray}
where $h_3$ is as in (\ref{h1}) and
$$
\widehat p^+(\al;x,y):={\rm e}^{-\lambda\,\al}\,p^+(\al;x,y),\quad
\widehat \eta_{x}(a,\al):={\rm e}^{-\lambda\,\al}\,\eta_{x}(a,\al).
$$
Because
$$
\widehat \eta_{x}(a,\al)=\P_x(H_a\in d\al,\, H_a<T)/d\al
$$
it is seen from (\ref{F}) that $F$ describes the finite dimensional distributions of
$X$ conditioned to hit $a$ at time $s$ before $T.$ For the
claim concerning the post-process we remark first that (\ref{G}) gives finite
dimensional distributions of the announced $h$-transform started from $a$ since
$$
\frac {\widehat \eta_{y_1}(a,v-s+v_1)}{\sqrt{2\lambda}}
= \lim_{x\downarrow a}\frac{\widehat p^+(v-s+v_1;x,y_1)}{h_3(x-a)}.
$$
Next notice that proceeding as above we can also compute the conditional probabilites for
$A_n$ and $B_{m}\circ\theta_v$ separately and deduce
\begin{eqnarray}
\label{b3}
&&\hskip-1cm
\P_0(A_n,  B_{m}\circ\theta_v\,|\,  I_T=a, H_{I}=s)\\
&&\hskip1cm
\nonumber
=\P_0(A_n\,|\,  I_T=a, H_{I}=s)\,\P_0(B_{m}\circ\theta_v\,|\,  I_T=a, H_{I}=s).
\end{eqnarray}
As is seen from (\ref{G}) the quantity 
$$
\P_0(B_{m}\circ\theta_v\,|\,  I_T=a, H_{I}=s)
$$
is a function of the difference $v-s$ only, and we find the desired
description of the post-process by letting $v\downarrow s$ (and
applying the
Lebesgue dominated convergence theorem). 
To remove the conditioning with respect to $H_I$ in (\ref{b3})  observe from
(\ref{levy3}) that
$$
\P_0(H_{I}\in ds\,|\, I_T=a)=\frac{ \widehat\eta_0(a,s)}{{\rm e}^{\,a\,\sqrt{2\lambda}}}\ ds
$$
and, hence,
\begin{eqnarray*}
&&
\hskip-1cm
\nonumber
\P_0(A_n\,|\,  I_T=a)=\widehat p^+(u_1;0,x_1)\,2dx_1 \\
&&\hskip1cm
\times\dots\cdot
\widehat p^+(u_n-u_{n-1};x_{n-1},x_n)\, 2dx_n\
\frac{{\rm e}^{-(x_n-a)\,\sqrt{2\lambda}}}{{\rm e}^{-a\,\sqrt{2\lambda}}},
\end{eqnarray*}
which means that the pre-process is as stated, and, moreover, the claimed conditional independence
holds. 

It is possible to prove claim 1 also via direct computations with finite dimensional distributions; however, we do not present this proof since, as seen below, the result is 
in the core of the approach with excursions. 
\hfill\break\hfill
{\sl c)} {\sl Excursion theoretical proof.}  The excursion process associated with the excursions
above the running minimum is defined similarly as the corresponding
process with running maximum in Section 2. Indeed, let 
for $a<0$ 
$$
H_{a-}:=\inf\{t\,:\,X_t<a\},
$$
and, if $H_{a-}>H_{a},$ 
$$
\Xi^-_a(u):=X_{u+H_{a}}-a,\quad  0\leq u<H_{a-}
  -H_{a}.
$$
Then the process 
$$
\Xi^-=\{ \Xi^-_a\,:\, a\leq 0\}
$$
is a homogeneous Poisson point process, and is called the excursion process for excursions
above the running minimum. We remark that $\Xi^-$ is identical in law
with the excursion process for
excursions from 0 of a reflecting Brownian motion. The Ito excursion measure associated with $\Xi^-$
is denoted by $n^-$ (for different descriptions of $n^-,$ see Revuz and
Yor \cite{revuzyor01}). 

Let $F_1$ and $F_2$ be measurable mappings
from $\cC(\R_+,\R)$ to $\R_+.$ Now we can write
\begin{eqnarray*}
&&
\hskip-.5cm
\Delta:=\E\left(F_1(X_u\,:\,u\leq H_I)\,F_2(X_{H_I+u}-I_T\,:\,u\leq
T-H_I)\right)\\
&&\hskip.5cm
=
\E\left(\sum_{a<0}F_1(X_u\,:\,u\leq H_a)\ F_2(\Xi^-_a(u)\,:\,u\leq
T-H_{a})\, {\bf 1}_{\{ H_a\leq T<H_{a-}\}} \right),
\end{eqnarray*}
where the sum is over all points of  $\Xi^-$ (but simplifies, for every
$\omega$ a.s., only to one term). Let $\cE$ denote the excursion space
and $\varepsilon$ a generic excursion.  
By the compensation formula for
Poisson point processes (see Bertoin \cite{bertoin96}) 
\begin{eqnarray}
\label{D11}
&&
\nonumber
\Delta=
\int_{-\infty}^0da\  
\E_0\Big(F_1(X_u\,:\,u\leq H_a)\,
{\bf 1}_{\{H_{a}<T\}}\\
&&
\nonumber
\hskip2cm
\times\int_{\cE}
F_2(\varepsilon_u\,:\, u\leq T-H_{a})\,
{\bf 1}_{\{T-H_{a}<\zeta\}}(\varepsilon)
\ n^-(d\varepsilon)\Big)\\
&&\hskip.5cm
=\int_{-\infty}^0da\  
\E_0\left(F_1(X_u\,:\,u\leq H_a)\,{\rm e}^{-\lambda H_a}\right)
\\
&&
\nonumber
\hskip2cm
\times
\int_{\cE} n^-(d\varepsilon)\ \E_0\left(
F_2(\varepsilon_u\,:\, u\leq T) 
\,
{\bf 1}_{\{T<\zeta\}}(\varepsilon)\right)
,
\end{eqnarray}
where the notation $\zeta(\varepsilon)$ is for the
life time of $\varepsilon$ and in the second step the fact that $T$ is an exponentially
distributed random variable independent of $X$ is used. Notice that 
(\ref{D11}) yields, when choosing $F_1\equiv 1,$ 
\begin{eqnarray*}
&&
\E_0\left(F_2(X_{H_I+u}-I_T\,:\,u\leq
T-H_I)\right)\\
&&
\hskip2cm
=\frac{1}{\sqrt{2\lambda}}\
\int_{\cE} n^-(d\varepsilon)\ \E_0\left(
F_2(\varepsilon_u\,:\, u\leq T) 
\,
{\bf 1}_{\{T<\zeta\}}(\varepsilon)\right).
\end{eqnarray*}
By absolute continuity, 
$$
\E_0\left(F_1(X_u\,:\,u\leq H_a)\,{\rm e}^{-\lambda H_a}\right)
=
\E_0^{-\sqrt{2\lambda}}\left(F_1(X_u\,:\,u\leq H_a)\right) \,{\rm
  e}^{a\,\sqrt{2\lambda}}
$$
and, since $-I_T$ is exponentially distributed with parameter
$\sqrt{2\lambda}$, we have 
\begin{eqnarray}
\label{b51}
&&\hskip-1cm
\nonumber
\Delta
=\int_{-\infty}^0da\  
\E_0^{-\sqrt{2\lambda}}\left(F_1(X_u\,:\,u\leq H_a)\right)\ \P(I_T\in da)
\\
&&
\hskip1.5cm
\times\,\frac{1}{\sqrt{2\lambda}}\ 
\int_{\cE} n^-(d\varepsilon)\ \E_0\left(
F_2(\varepsilon_u\,:\, u\leq T) 
\,
{\bf 1}_{\{T<\zeta\}}(\varepsilon)\right)
.
\end{eqnarray}
Consequently, the processes $\{X_u\,:\,u\leq H_I\}$ and
$\{X_{H_I+u}-I_T\,:\,u\leq T-H_I\}$ are independent and, hence, also 
$\{X_u\,:\,u\leq H_I\}$ and
$\{X_{T-t}-X_T\,:\,u\leq T-H_I\}$ are independent. Moreover,
$\{X_u\,:\,u\leq H_I\}$ given $I_T=a$ is identical in law with
BM($-\sqrt{2\lambda}$) killed at the first hitting time of $a.$
To prove that $\{X_{T-t}-X_T\,:\,u\leq T-H_I\}$ given $X_T-I_T=b$ is identical in law with
BM($-\sqrt{2\lambda}$) killed at the first hitting time of $-b$ observe first that 
\begin{eqnarray*}
&&
\int_{\cE} n^-(d\varepsilon)\ \E_0\left(
F_2(\varepsilon_u\,:\, u\leq T) 
\,
{\bf 1}_{\{T<\zeta\}}(\varepsilon)\right)\\
&&\hskip3cm
= \lambda\,
\int_0^\infty dt\ {\rm e}^{-\lambda t}\ 
n^-(F_2(\varepsilon_u\,:\, u\leq t)
\,
{\bf 1}_{\{t<\zeta\}}(\varepsilon)).
\end{eqnarray*}
Next we claim that 
\begin{eqnarray}
\label{b52}
&&\nonumber
\int_0^\infty dt\, {\rm e}^{-\lambda t}\,n^-(F_2(\varepsilon_u\,:\,
u\leq t)
\,
{\bf 1}_{\{t<\zeta\}}(\varepsilon)).
\\
&&\hskip2.5cm
=
2\,\int_0^\infty db\ \E_b\left({\rm e}^{-\lambda
  H_0}\,F_2(X_{H_0-u}\,:\, u\leq H_0)\right).
\end{eqnarray}
Indeed, (\ref{b52}) for $\lambda=0$ is formula 5 in 
Biane and Yor \cite{bianeyor87} Th\'eor\`eme 6.1 p. 79 
and the validity for $\lambda>0$ is easily verified by inspecting 
the proof in \cite{bianeyor87} p. 79. Hence, by spatial
symmetry,
\begin{eqnarray}
\label{b525}
&&\nonumber
\E_0\left(F_2(X_{H_I+u}-I_T\,:\,u\leq
T-H_I)\right)\\
&&\hskip2cm
=
\sqrt{2\lambda}
\,\int_0^\infty db\ \E_b\left({\rm e}^{-\lambda H_0}\,F_2(X_{H_0-u}\,:\, u\leq H_0)\right)\\
&&\hskip2cm\nonumber
=
\sqrt{2\lambda}
\,\int_0^\infty db\ \E_0\left({\rm e}^{-\lambda H_{-b}}\,F_2(b+X_{H_{-b}-u}\,:\, u\leq H_{-b})\right).
\end{eqnarray}
Reversing here time and using absolute continuity yield 
\begin{eqnarray*}
&&
\E_0\left(F_2(X_{T}-X_{T-u}\,:\,u\leq
T-H_I)\right)\\
&&\hskip2cm
=
\sqrt{2\lambda}
\,\int_0^\infty db\ \E_0\left({\rm e}^{-\lambda
  H_{b}}\,F_2(X_s\,:\, s\leq H_{b})\right)\\
&&\hskip2cm
=
\sqrt{2\lambda}
\,\int_0^\infty db\ \E_0^{\sqrt{2\lambda}}\left(F_2(X_s\,:\, s\leq H_{b})\right)
{\rm e}^{-b\,\sqrt{2\lambda}}
\end{eqnarray*}
which proves the first claim of the theorem. 

It remains to verify claim 2 (iii). For a fixed $t>0$ we obtain from (\ref{b525}) 
\begin{eqnarray*}
&&\nonumber
\widehat\Delta:=\E_0\left(F_2(X_{H_I+u}-I_T\,:\,u\leq t)\,{\bf 1}_{\{t\leq T-H_I\}}\right)\\
&&\hskip.5cm
=
\sqrt{2\lambda}
\,\int_0^\infty db\ \E_b\left({\rm e}^{-\lambda H_0}\,F_2(X_{H_0-u}\,:\, u\leq t)
\,{\bf 1}_{\{t\leq H_0\}} \right).
\end{eqnarray*} 
According to Williams' time reversal theorem the process $\{X_{H_0-u}:\, 0\leq u< H_0\}$
under $\P_b$ is identical in law with BES(3) started from 0 and killed at the last 
exit time at $b.$ Consequently, letting $\gamma_b$ denote the last exit time       
we have 
\begin{eqnarray*}
&&\nonumber
\widehat\Delta=
\sqrt{2\lambda}
\,\int_0^\infty db\ \Q_0\left({\rm e}^{-\lambda \gamma_b}\,F_2(X_{u}\,:\, u\leq t)
\,{\bf 1}_{\{t\leq \gamma_b\}} \right)\\
&&\hskip.5cm
=
\sqrt{2\lambda}\,{\rm e}^{-\lambda\, t}
\,\int_0^\infty db\ \Q_0\Big(F_2(X_{u}\,:\, u\leq t)\,\Q_{X_t}\left(
{\rm e}^{-\lambda \gamma_b} \,{\bf 1}_{\{\gamma_b>0\}}\right) \Big)
\end{eqnarray*} 
by the Markov property. The distribution of $\gamma_b$ is well known (see Pitman and Yor 
\cite{pitmanyor81}) and  it holds
$$
\Q_{r}\left(
{\rm e}^{-\lambda \gamma_b} \,{\bf 1}_{\{\gamma_b>0\}}\right)= \frac 1{r\sqrt{2\lambda}}\ 
\sh((b\wedge r)\sqrt{2\lambda})\ {\rm e}^{-(b\vee r)\,\sqrt{2\lambda}}.
$$
We have now
\begin{eqnarray*}
&&\nonumber
\widehat\Delta=
\sqrt{2\lambda}\,{\rm e}^{-\lambda\, t}
\,\Q_0\left(F_2(X_{u}\,:\, u\leq t)\,\int_0^\infty db\ \Q_{X_t}\left(
{\rm e}^{-\lambda \gamma_b} \,{\bf 1}_{\{\gamma_b>0\}}\right) \right)
\\&&\hskip.5cm
=
{\rm e}^{-\lambda\, t}
\,\Q_0\left(F_2(X_{u}\,:\, u\leq t)\,\frac{1-{\rm e}^{-X_t \sqrt{2\lambda}}}
{X_t\,\sqrt{2\lambda}}\right)\\
&&\hskip.5cm
=
\Q_0\left(F_2(X_{u}\,:\, u\leq t)\ \frac{h_3(X_t)}
{X_t\,\sqrt{2\lambda}}\ {\bf 1}_{\{t\leq T\}}\right),
\end{eqnarray*} 
where $h_3$ is as in (\ref{h1}). Define 
$h_4(x):=h_3(x)/x$ for $x>0$ and $h_4(0):=\sqrt{2\lambda},$
and notice that $h_4$ is right continuous at 0. Consequently, 
\begin{eqnarray}
\label{bes0}
&&\nonumber
\E_0\left(F_2(X_{H_I+u}-I_T\,:\,u\leq t)\,{\bf 1}_{\{t\leq T-H_I\}}\right)
\\&&\hskip3cm
=
\Q_0\left(F_2(X_{u}\,:\, u\leq t)\,\frac{h_4(X_t)}{h_4(0)}\ {\bf 1}_{\{t\leq T\}}\right).
\end{eqnarray} 
It is easy to verify that for $x>0$
$$
\cG^{R}h_4(x)-\lambda h_4(x)=-\frac\lambda x<0,
$$
where $\cG^{R}$ is the infinitesimal generator of BES(3), 
see (\ref{e131}). Using It\^o's formula it is seen  that $h_4$ is $\lambda$-excessive for BES(3). 
Consequently, the post $H_I$-process is identical in law with the Doob $h$-transform with 
$h=h_4$ of exponentially killed BES(3). The generator of this transform can be computed in the usual way, and is seen to coincide when $a=0$ with 
$\cG^Z$ given in (\ref{h0}). 
\end{proof}

\begin{remark} 
\label{rem311}
{\bf 1.} {\sl Informally, given $I_T=a$ the post-$H_I$-process is
identical in law with a 
Brownian motion killed at time $\zeta:=T\wedge H_a$ and conditioned by the event $X_{\zeta-}>a.$ 
Let $\P^Z$ denote the probability measure associated with the diffusion $Z$ introduced in Theorem \ref{deco1}. Then for any finite stopping time $U$ and $\Lambda_U\in \cF_U$ 
it holds 
\begin{eqnarray*}
&&
\P_x^{Z}(\Lambda_U) =\frac 1 {h_3(x-a)}\E_x\left(h_3(X_U-a); \Lambda_U, U<T\wedge H_a\right)\\
&&\hskip1.6cm
=\P_x\left(\Lambda_U, U<T\wedge H_a\,|\, X_{\zeta-}> a\right).
\end{eqnarray*}
The diffusion $Z$ can alternatively, as is seen at the end of the above proof, be described as 
the Doob $h$-transform with $h=h_4$ of exponentially killed BES(3). 
}

\hfill\break\hfill
{\bf 2.} {\sl From Theorem \ref{deco1} it is clear that
$\{X_{t}:\, 0\leq t< H_I\}$ and  $\{X_{H_I+t}-I_T:\, 0\leq t< T-H_I\}$
are independent. In particular, $X_T-I_T$ and $I_T$ are independent.
This last property is also easily verified by analyzing the joint distribution of
$X_T$ and $I_T$ 
$$ 
\P_0(I_T\in da, X_T\in dz)=2\lambda\, {\rm e}^{2a\,\sqrt{2\lambda}}
\, {\rm e}^{-z\,\sqrt{2\lambda}}\, da\,dz,\quad a<0,\ a<z,
$$
which is obtained from (\ref{levy2}) by integrating
with respect to $s.$ Moreover, $X_T-I_T$ and $-I_T$ are seen to be identically
distributed the common distribution being the exponential distribution
with parameter $\sqrt{2\lambda}$ (cf. also (\ref{I1}) below). 
}\hfill\break\hfill
{\bf 3.} {\sl Notice also the fact that $H_I$ and $T-H_I$ are
independent and identically Gamma$(1/2,\lambda)$-distributed, i.e.,
$$
\P_0(H_I\in du)=\P_0(T-H_I\in
du)=\frac{\sqrt{\lambda}}{\sqrt{\pi u}}
\, {\rm e}^{-\lambda\, u}\, du,\quad  u>0.
$$
}
\end{remark}

Next we recall the formulas (see \cite{borodinsalminen02} p. 173)
\begin{equation}
\label{joint1}
\P_0(I_T\in da,S_T\in db)=\lambda\ \frac{\ch(\frac 12
  (b+a)\sqrt{2\lambda})}
{\ch^3(\frac 12 (b-a)\sqrt{2\lambda})}\ da\,db,\quad  a<0<b, 
\end{equation}
and 
\begin{equation}
\label{I1}
\P_0(I_T\in da)=\sqrt{2\lambda}\,\e^{a\sqrt{2\lambda}}\, da\quad  a<0.
\end{equation}
We need also the following result derived as a corollary to Theorem
\ref{deco1} but which can also be deduced from the
formula 1.1.28.2 p. 191 in \cite{borodinsalminen02}.
\begin{corollary}
\label{hit1}
{\sl For $a<0<b$ 
\begin{eqnarray}
\label{hit2} 
&&\hskip-1cm
\nonumber
\P_0(H_I<H_S,I_T\in da,S_T\in db)\\
&&\hskip1.8cm
=
2\lambda\ \frac{\left(\ch(
  (b-a)\sqrt{2\lambda})-1\right)\,\sh(b\sqrt{2\lambda})}
{\sh^3((b-a)\sqrt{2\lambda})}\ da\,db.
\end{eqnarray} 
}
\end{corollary}
\begin{proof}
The conditional independence stated in Theorem \ref{deco1} yields 
\begin{eqnarray}
\label{HH}
&&\hskip-1cm
\nonumber
\P_0(H_I<H_S,S_T\in db\,|\, I_T=a)\\
&&\hskip1.8cm
\nonumber
=
\P_0(H_I<H_S,S_{H_I,T}\in db\,|\, I_T=a)\\
&&\hskip1.8cm
=
\P_0(S_{H_I,T}\in db\,|\, I_T=a)\ 
\P_0(S_{H_I}<b\,|\, I_T=a),
\end{eqnarray} 
where
$$
 S_{H_I,T}:=\sup\{X_t:\, H_I\leq t\leq T\}.
$$
From the description of the pre-$H_I$ process we have
using the well known formula, see e.g. \cite{borodinsalminen02} p. 309,
\begin{eqnarray*} 
&&
\P_0(S_{H_I}<b\,|\, I_T=a)=\P_0^{-\sqrt{2\lambda}}(H_a<H_b)\\
&&\hskip3.8cm
=\e^{-a\sqrt{2\lambda}}\ \frac{\sh (b\,\sqrt{2\lambda})}
{\sh ((b-a)\,\sqrt{2\lambda})}.
\end{eqnarray*} 
Using the $h$-transform description of $Z$ (cf. Remark \ref{rem311}) 
\begin{eqnarray}
\label{h/h}
&&\hskip-2.2cm
\nonumber
\P_0(S_{H_I,T}>b\,|\, I_T=a)= 
\P_a^{Z}(H_b<\infty)\\
&&\hskip2cm
\nonumber
=
\lim_{x\downarrow a}\frac{h_3(b-a)}{h_3(x-a)}\,\P_x(H_b<H_a,H_b<T)\\
&&\hskip2cm
\nonumber
=
\lim_{x\downarrow  a}
\frac{h_3(b-a)}{h_3(x-a)}\,\frac{\sh((x-a)\sqrt{2\lambda})}
{\sh((b-a)\sqrt{2\lambda})}\\
&&\hskip2cm
=
\frac{1-\e^{-(b-a)\sqrt{2\lambda}}}
{\sh((b-a)\sqrt{2\lambda})}.
\end{eqnarray} 
The formula (\ref{hit2}) is now obtained from (\ref{HH}) when
multiplying with the density of $I_T$ given in (\ref{I1}).
\end{proof}

We proceed by refining the path decomposition presented in Theorem
\ref{deco1}. Recall that $\Q^{\mu}_{\,x}$ with $\mu, x\geq 0$ denotes 
the measure under which the coordinate 
process $X=\{X_t\,:\, t\geq 0\}$ is a BES(3,$\mu$) started from $x$ (see 
(\ref{e13}) for the generator). 

\begin{theorem}
\label{deco2}
{\sl  Conditionally on $H_I<H_S,\
I_T=a,\ S_T=b$ it holds that 
\begin{enumerate}
\item the pre-$H_I$-process $\{X_t:\, 0\leq t\leq H_I\}$ is
  identical in law with $\{b-X_t:\, 0\leq t\leq
  H_{b-a}\}$ under $\Q^{\sqrt{2\lambda}}_{\,b}.$
\item the intermediate process $\{X_{H_I+t}:\, 0\leq t\leq H_S-H_I\}$
is identical in law with $\{a+ X_t\,:\, 0\leq t\leq
H_{b-a}\}$ under $\Q^{\sqrt{2\lambda}}_{\,0}.$
\item
the post-$H_S$-process $\{X_{H_S+t}:\, 0\leq t\leq T-H_S\}$
is identical in law with a diffusion $Z^\circ$ started from $b$ and having the generator
$$
{\cal G}^{Z^\circ}u(x)=\frac 12 u''(x)
+\frac{g'(x)}{g(x)}u'(x)-\frac{\lambda g(x)-\frac 12\,g''(x)}{g(x)}u(x),
$$
where $a<x<b$ and
\begin{equation}
\label{g1}
g(x)=\sh((b-a)\sqrt{2\lambda})-\sh((b-x)\sqrt{2\lambda})-\sh((x-a)\sqrt{2\lambda}),
\end{equation}
\item  the  pre-$H_I$-process, the intermediate process and  the post-$H_S$-process are
independent.
\end{enumerate}
}
\end{theorem}

\begin{proof} From Theorem \ref{deco1}, using symmetry of BM, 
it is seen that conditionally on $S_T=b$
\begin{description}
\item{P1.} the processes  $\{X_t:\, 0\leq t\leq H_{S}\}$ and 
 $\{X_{H_{S}+t} :\, 0\leq t\leq T- H_S\}$ are independent, 
\item{P2.} the pre-$H_S$-process $\{X_{t} :\,
  0\leq t\leq H_S\}$ is distributed as  BM($\sqrt{2\lambda}$) started from 0 and 
killed when it reaches $b,$
\item{P3.} the post-$H_S$-process $\{X_{H_S+t} :\,
  0\leq t\leq T-H_S\}$ is a diffusion $Z^{\downarrow}$ started from $b$ and having the generator
\begin{equation}
\label{h00}
{\cal G}^{Z^\downarrow}u(x)=\frac 12 u''(x) +\frac{h'(b-x)}{h(b-x)}u'(x)-\frac{\lambda}{h(b-x)}u(x).
\end{equation}  
\end{description}

For the pre-$H_S$-process we apply the path decomposition theorem of Brownian
motion with drift as given in Tanr\'e and Vallois \cite{tanrevallois04} Proposition 3.2.
To recall this, let 
$$
I_{H_b}:=\inf\{X_t:\, 0\leq t\leq H_b\}
\quad
{\rm and}
\quad
\widehat H_{b}:=\inf\{t:\,  X_t = I_{H_b}\}.
$$
Then, under $\P^\mu_0,$ conditionally on $I_{H_b}=a$
\begin{description}
\item{P4.} the processes  $\{X_t:\, 0\leq t\leq \widehat H_{b}\}$ and 
 $\{X_{\widehat H_{b}+t} :\, 0\leq t\leq H_b-\widehat  H_b\}$ are independent, 
\item{P5.}   the pre-$\widehat H_{b}$-process $\{X_{t} :\, 0\leq t\leq \widehat H_{b}\}$
is identical in law with  BM($-\mu$) started from 0, conditioned to
hit $a$ before $b,$ and killed when it hits $a.$
\item{P6.} the post-$\widehat H_{b}$-process $\{X_{\widehat H_{b}+t} :\, 0\leq t\leq
  H_b-\widehat H_b\}$ is identical in law with $\{a+X_t:\, 0\leq t\leq \widehat H_{b-a}\}$
under $\Q^\mu_0.$
\end{description}

Let $F_1,$ $F_2,$ and $F_3$ be measurable mappings
from $\cC(\R_+,\R)$ to $\R_+,$ and introduce   
$$
F_1:=F_1(\{X_{t} :\,
  0\leq t\leq H_I\}),
\quad
F_2:=F_2(\{X_{H_I+t} :\,
  0\leq t\leq H_S-H_I\}),
$$
and
$$
F_3:=F_3(\{X_{H_S+t} :\,
  0\leq t\leq T-H_S\}),
$$
where it is assumed that $H_I<H_S.$ Notice that
\begin{equation}
\label{equiv1}
H_I\leq H_S\quad \Leftrightarrow\quad  I_{H_S,T}:=\inf\{X_u:\, H_S\leq
u\leq T\}> I_{H_S}.
\end{equation}
and 
$$
 I_{H_S,T}> I_{H_S}\quad \Leftrightarrow\quad  I_{H_S}= I_T
$$ 
Consider now
\begin{eqnarray*}
&&
\E_0\left(F_1F_2F_3\ ;\  I_T\in da, H_I<H_S \,|\,S_T=b \right)
\\
&&\hskip2cm
=
\E_0\left(F_1F_2F_3\ ;\ I_{H_S}\in da, I_{H_S,T}>a \,|\,S_T=b \right)
\\
&&\hskip2cm
=
\E_0\left(F_1F_2\ ;\ I_{H_S}\in da\,|\,S_T=b \right)
\E_0\left(F_3\ ;\ I_{H_S,T}>a \,|\,S_T=b \right),
\end{eqnarray*}
where in the second step the conditional independence (see P1) is applied. 
The first term on the right hand side can be analyzed via properties P4, P5 and P6. In particular, claims 1 and 2 of Theorem \ref{deco2} follow with the help of Lemma \ref{lem1}. Moreover, P1 and P4 yield claim 4. 

Next we prove claim 3. For this 
let $\P^\downarrow_x$ denote the measure in the canonical setting associated with $Z^{\downarrow}$ when started from $x<b.$
Hence,
\begin{eqnarray*}
&&
\E_0\left(F_3(\{X_{H_S+u} :\,
  0\leq u\leq  t\})\,{\bf 1}_{\{t\leq T-H_S\}}\, |\, I_{H_S,T}>a ,S_T=b \right)
\\
&&\hskip2cm
=
\E^\downarrow_b\left(F_3(\{X_u:\,
  0\leq u\leq  t\})\,{\bf 1}_{\{t\leq \zeta\}}\, |\, I_\zeta>a  \right),
\end{eqnarray*}
where $\zeta$ denotes the life time and $I_\zeta$ the global infimum. By the Markov property,
\begin{eqnarray*}
&&
\E^\downarrow_b\left(F_3(\{X_u:\,
  0\leq u\leq  t\})\,;\,t\leq \zeta\,,\, I_\zeta>a  \right)\\
&&\hskip1cm
=
\E^\downarrow_b\left(F_3(\{X_u:\,
  0\leq u\leq  t\})\,{\bf 1}_{\{t<H_a\}}\,\P^{\downarrow}_{X_t
}\left( I_\zeta>a\right)  \right).
\end{eqnarray*}
Clearly, for $x\leq b$ 
$$
\P^{\downarrow}_x\left( I_\zeta>a\right)=\P^{\downarrow}_x\left( H_a=\infty\right)=
1-\P^{\downarrow}_x\left( H_a<\infty\right),
$$
and computing as in (\ref{h/h}) we obtain
$$
\P^{\downarrow}_x\left( H_a<\infty\right)=
\frac{h_3(b-a)}{h_3(b-x)}\,\frac{\sh((b-x)\sqrt{2\lambda})}
{\sh((b-a)\sqrt{2\lambda})}
$$
with $h_3$ as in (\ref{h1}). Consequently, after some elementary manipulations,
for $x\in(a,b)$ 
$$
\P^{\downarrow}_x\left( I_\zeta>a\right)
=h_5(x):=
\frac {g(x)}{h_3(b-x)\,\sh((b-a)\sqrt{2\lambda})}
$$
with $g$ as in (\ref{g1}). Since  $Z^\downarrow$ is the Doob $h$-transform with $h=h_3(b-\cdot)$ of BM killed at time $T\wedge H_b$ it follows that
\begin{eqnarray*}
&&
\E^\downarrow_b\left(F_3(\{X_u:\,
  0\leq u\leq  t\})\,;\,t\leq \zeta\,|\, I_\zeta>a  \right)\\
&&\hskip1cm
=
\E^\downarrow_b\left(\frac{h_5(X_t)}{h_5(b)}\,F_3(\{X_u:\,
  0\leq u\leq  t\})\,;\,t\leq \zeta \right)
\\
&&\hskip1cm
=
\lim_{x\uparrow b}\E_x\left(\frac{h_3(b-X_t)\,h_5(X_t)}{h_3(b-x)\,h_5(x)}\,F_3(\{X_u:\,
  0\leq u\leq  t\})\,;\,t\leq T\wedge H_b\wedge H_a \right)
\\
&&\hskip1cm
=
\lim_{x\uparrow b}\E_x\left(\frac{g(X_t)}{g(x)}\,F_3(\{X_u:\,
  0\leq u\leq  t\})\,;\,t\leq T\wedge H_b\wedge H_a \right),
\end{eqnarray*}
and this yields the description 3 of the post process. 
\end{proof}

\begin{remark}
\label{info}
{\bf 1.} Informally, the post-$H_S$ process is identical in law with a
Brownian motion killed at time $\zeta:=T\wedge H_a\wedge H_b$ and conditioned by the
event $X_{\zeta-}\in (a,b).$ Indeed, from \cite{borodinsalminen02} 3.0.1
p. 212 
\begin{eqnarray*}
&&
\P_x\left(X_\zeta\in(a,b)\right)=\P_x\left(T<H_a\wedge H_b\right)
\\
&&\hskip3cm
= \E_x\left(1-\exp(-\lambda H_a\wedge H_b) \right)
\\
&&\hskip3cm
= 1-\frac{\sh((b-x)\sqrt{2\lambda})-\sh((x-a)\sqrt{2\lambda})}{ \sh((b-a)\sqrt{2\lambda})}
\\
&&\hskip3cm
=\frac{g(x)}{ \sh((b-a)\sqrt{2\lambda})}.
\end{eqnarray*}
\hfill\break\hfill
{\bf 2.} It can be proved that the process  $\{X_{T-t}:\, 0\leq t\leq T-H_S\}$
, i.e., 
the time reversal of the post-$H_S$-process is identical in law with 
$\{a+X_{t}:\, 0\leq t\leq H_{b-a}\}$  under $\Q^{\sqrt{2\lambda}}_{\,b-a-\xi}$
with $\xi$ a random variable independent of $X$ having the density 
$$
f_\xi(x)=\sqrt{2\la}\
\e^{-x\,\sqrt{2\la}}/(1-\e^{-(b-a)\,\sqrt{2\la}}),\quad 0<x<b-a.
$$
Consequnetly, for  $0<x<b-a$ 
$$
\P(S_T-X_T\in dx\ |\ S_T=b, I_T=a, H_I<H_S)=f_\xi(x)\,dx.
$$
\end{remark}

 \section{Maximum increase and decrease}
\label{sec3}
In this section we apply the results from the previous sections to find
an expression for the joint distribution of the maximum increase and
the maximum decrease of Brownian motion up to an independent
exponential time $T$. As stated in the Introduction the law of 
$(D^+_t,D^-_t)$ under $P_x$ does not depend on $x$ and, hence, we write in the sequel $\P$ instead of 
$\P_x.$ 

For $s<t$ define 
the maximum increase on the interval $(s,t)$ via 
$$
D^+_{s,t}:=\sup_{s\leq u\leq v\leq t}\left(X_v -X_u\right)
$$
and, analogously, the maximum decrease
$$
D^-_{s,t}:=\sup_{s\leq u\leq v\leq t}\left(X_u -X_v\right).
$$
With these new notations we have $D^-_{T}=D^-_{0,T}$ and $D^+_{T}=D^+_{0,T}.$ 
 
The path decomposition given in Theorem \ref{deco2} leads us to the
following 
\begin{proposition}
\label{prop41}
{\bf 1.}
For $-a<d<b-a$
\begin{eqnarray*}
&&\hskip-1cm
\P(D^-_{H_I}<d\ |\ H_I<H_S, I_T=a,S_T=b)
\\
&&\hskip3cm
=
\frac{S^{\nu}(a+d)(S^{\nu}(b)- S^{\nu}(a))}{S^{\nu}(b)(S^{\nu}(a+d)-S^{\nu}(a))}
\\
&&\hskip3cm
= \frac{\sh((a+d)\sqrt{2\la})\,\sh((b-a)\sqrt{2\la})}{\sh(d\,\sqrt{2\la})\,\sh(b\,\sqrt{2\la})}
\\
&&\hskip3cm
=: f_1(d;a,b),
\end{eqnarray*}
where $S^{\nu}$ is as (\ref{e110}) with $\nu=-\sqrt{2\lambda},$ i.e.,
$$
S^{\nu}(x):=\frac{1}{2\sqrt{2\lambda}}({\rm e}^{2\sqrt{2\lambda}
  x}-1)
=
\frac{1}{\sqrt{2\lambda}}\ {\rm e}^{\sqrt{2\lambda}\,x}\ \sh(x\,\sqrt{2\lambda}).
$$
\hfill\break\hfill
{\bf 2.}
For $0<d<b-a$
\begin{eqnarray*}
&&
\P(D^-_{H_I,H_S}<d\ |\ H_I<H_S, I_T=a,S_T=b)\\
&&\hskip3cm
=\disp{\frac{S^{\nu}(b-a)}{S^{\nu}(d)}}
\exp\left(-\disp{\frac{b-a-d}{S^{\nu}(d)}}-2\sqrt{2\lambda}(b-a-d)\right)\\
&&\hskip3cm
=
\frac{\sh((b-a)\sqrt{2\la})}{\sh(d\,\sqrt{2\la})}\exp\left(-(b-a-d)\sqrt{2\la}\ 
 \cth(d\,\sqrt{2\la})\right)
\\
&&\hskip3cm
=: f_2(d;a,b)
\end{eqnarray*}
\hfill\break\hfill
{\bf 3.}
For $0<d<b-a$
\begin{eqnarray*}
&&
\P(D^-_{H_S,T}<d\ |\ H_I<H_S, I_T=a,S_T=b)
\\
&&\hskip3cm
=\frac{\sh((b-a)\sqrt{2\la})(\ch(d\sqrt{2\la})-1)}{\sh(d\sqrt{2\la})(\ch((b-a)\sqrt{2\la})-1)}.
\\
&&\hskip3cm
=: f_3(d;a,b)
\end{eqnarray*}
Moreover, conditionally on $I_T,$ $S_T$ and $H_I<H_S$ the variables
$D^-_{H_I},$ $D^-_{H_I,H_S},$ and $D^-_{H_S,T}$ are independent.   
\end{proposition}
\begin{proof} 
{\sl Claim 1.}  Because
$$
D^-_{T}=\sup_{0\leq v\leq t}\left(\sup_{0\leq u\leq v}X_u -X_v\right)
$$
and since, conditionally on $I_T=a,$ $H_I$ is the first hitting time
of $a$ we have 
 $$
D^-_{H_I}=\sup_{0\leq v\leq H_I} X_u -a 
$$
and, hence, for $-a<d<b-a$
$$
\{D^-_{H_I}<d\}=\{H_a<H_{a+d}\}.
$$
From  Theorem \ref{deco2} and Lemma \ref{lem1} it follows that 
\begin{eqnarray*}
&&
\P(D^-_{H_I}<d\ |\ H_I<H_S, I_T=a,S_T=b)
\\
&&\hskip2.5cm
=\Q^{\sqrt{2\lambda}}_b(H_{b-a}<H_{b-a-d})
\\
&&\hskip2.5cm
=\P_b^{\sqrt{2\lambda}}(H_{b-a}<H_{b-a-d}\ |\
H_{b-a}<H_{0}).
\end{eqnarray*}
Since $\{b-X_t\,:\, t\geq 0\}$  under
$\P_b^{\sqrt{2\lambda}}$ is distributed as $\{X_t\,:\, t\geq 0\}$ under $\P_0^{-\sqrt{2\lambda}},$
we have 
\begin{eqnarray*}
&&
\P(D^-_{H_I}<d\ |\ H_I<H_S, I_T=a,S_T=b)
\\
&&\hskip2.5cm
=\P_0^{-\sqrt{2\lambda}}(H_{a}<H_{a+d}\ |\
H_{a}<H_{b})
\\
&&
\hskip2.5cm
=
\P_0^{-\sqrt{2\lambda}}(H_a<H_{a+d})/\P_0^{-\sqrt{2\lambda}}(
H_a<H_{b})
\\
&&
\hskip2.5cm
=
\frac{S^{\nu}(a+d)-S^{\nu}(0)}{S^{\nu}(a+d)-S^{\nu}(a)}\,
\frac{S^{\nu}(b)- S^{\nu}(a)}{S^{\nu}(b)-S^{\nu}(0)},
\end{eqnarray*}
where the fact that  $S^{\nu}$ is the scale function  of
BM($-\sqrt{2\lambda}$) is used. Observing that for $y>x$
$$
S^{\nu}(y)-S^{\nu}(x)
=
\frac{1}{\sqrt{2\lambda}}\ {\rm e}^{\sqrt{2\lambda}\,(y+x)}\
\sh((y-x)\,\sqrt{2\lambda})
$$
leads immediately to the claimed formula. 
\hfill\break\hfill
{\sl Claim 2.} This follows directly from Proposition \ref{prop2} and Theorem \ref{deco2}.
\hfill\break\hfill
{\sl Claim 1.} Using  again  Theorem \ref{deco2}
\begin{eqnarray*}
&&
\P(D^-_{H_S,T}<d\ |\ H_I<H_S, I_T=a,S_T=b)
=\P_b^{Z^\circ}(H_{b-d}=+\infty).
\end{eqnarray*}
By the $h$-transform description of $Z^\circ$ (cf. Remark \ref{info}), 
\begin{eqnarray*}
&&
\P_b^{Z^\circ}(H_{b-d}<+\infty)=\lim_{x\uparrow
  b}\P_x^{Z^\circ}(H_{b-d}<+\infty)
\\
&&\hskip3.4cm
=\lim_{x\uparrow
  b}\frac{1}{g(x)}\E_x(g(X_{H_{b-d}})\,;\, H_{b-d}<T\wedge H_b\wedge H_a)\\
&&
\hskip3.4cm
=\lim_{x\uparrow
  b}\frac{g(b-d)}{g(x)}\E_x(\e^{-\la\, H_{b-d}}\,;\, H_{b-d}<H_b)\\
&&
\hskip3.4cm
=\lim_{x\uparrow
  b}\frac{g(b-d)}{g(x)}\frac{\sh((b-x)\sqrt{2\lambda})}{\sh(d\sqrt{2\lambda})},
\end{eqnarray*}
where $g$ is as in (\ref{g1}), i.e., for $a<x<b$
\begin{equation}
\label{g11}
g(x)=\sh((b-a)\sqrt{2\lambda})-\sh((b-x)\sqrt{2\lambda})-\sh((x-a)\sqrt{2\lambda}).
\end{equation}
Straightforward computations yield now the claimed formula.
\end{proof}

Applying symmetry properties of Brownian motion gives us the following
formulas
\begin{equation}
\label{e41}
\P(D^+_{T}<\al\,,\,D^-_{T}<\beta)=
\P(D^+_{T}<\beta\,,\,D^-_{T}<\al),
\end{equation}
\begin{eqnarray}
\label{e4105}
&&
\P(D^+_{T}<\al\,,\,D^-_{T}<\beta\,,\, H_S<H_I)\\
&&\nonumber
\hskip3cm
=
\P(D^+_{T}<\beta\,,\,D^-_{T}<\al\,,\, H_I<H_S),
\end{eqnarray}
and 
\begin{eqnarray}
\label{e411}
&&\hskip-1cm
\P(D^+_{T}<\al\,,\,D^-_{T}<\beta)
=
\P(D^+_{T}<\al\,,\,D^-_{T}<\beta\,,\, H_I<H_S)\\
&&
\nonumber
 \hskip5cm 
+
\P(D^+_{T}<\beta\,,\,D^-_{T}<\al\,,\, H_I<H_S).
\end{eqnarray}
The  distribution and the  density function of
$(D^+_T,D^-_T)$ will result as a corollary of the next theorem 
when applying the formulas (\ref{e41}), (\ref{e4105}), and (\ref{e411}).

\begin{theorem}
\label{thm41}
For $\beta>0$  and $\varphi:\R_+\mapsto\R_+$ bounded and measurable it
holds 
\begin{eqnarray}
\label{en41}
&&
\hskip-.5cm
\E\left(\varphi(D^+_{T})\,;\, D^-_{T}<\beta\, ,\,H_I<H_S\right)\\
&&\nonumber
=
\sqrt{2\la} \frac{
  (\ch(\beta\sqrt{2\la})-1)^2}{\sh^3(\beta\sqrt{2\la})}
\int_\beta^\infty dx \,\varphi(x)\, 
 \exp\left(-(x-\beta)\sqrt{2\la}\ \cth(\beta\,\sqrt{2\la})\right)\\
&&\nonumber
\hskip2.5cm
+
\sqrt{2\la}\int_0^\beta dx\, \varphi(x)\,  \frac{
  (\ch(x\sqrt{2\la})-1)^2}{\sh^3(x\sqrt{2\la})}.
\end{eqnarray}
\end{theorem}

\begin{proof}
Notice first that if $H_I<H_S$ then 
$$
D^+_{T}=S_T-I_T\quad {\rm and}\quad  D^-_{T}=D^-_{H_I}\vee
D^-_{H_I,H_S} \vee D^-_{H_S,T},
$$
and, therefore, 
\begin{eqnarray*}
&&
\hskip-.8cm
\Delta:=
\E\left(\varphi(D^+_{T})\,;\, D^-_{T}<\beta\, ,\,H_I<H_S\right)\\
&&
\hskip.5cm
=\E\left(\varphi(S_T-I_T)\,;\, D^-_{H_I}<\beta\,,\,
D^-_{H_I,H_S}<\beta\,,\, D^-_{H_S,T}<\beta\,,\,H_I<H_S\right).
\end{eqnarray*}
Taking the conditional expectation with respect to the
$\sigma$-algebra generated by $(I_T,S_T,{\bf 1}_{\{H_I<H_S\}})$ and
using Proposition \ref{prop41} we obtain 
\begin{eqnarray}
\label{en43}
&&
\Delta=\int_{-\infty}^0 da\int_0^{+\infty}db\, \varphi(b-a)  \,
f(a,b)\\
&&\nonumber
\hskip2cm
\times {\bf 1}_{\{\beta>-a\}}\left( f_1(\beta;a,b)\, f_2(\beta;a,b)\, f_3(\beta;a,b)\,{\bf 1}_{\{\beta<b-a\}}
+{\bf 1}_{\{\beta>b-a\}}\right)
\end{eqnarray}
with $f_i,\, i=1,2,3,$ as in Proposition \ref{prop41} and 
\begin{eqnarray*}
&&
f(a,b)=2\lambda\ \frac{\left(\ch(
  (b-a)\sqrt{2\lambda})-1\right)\,\sh(b\sqrt{2\lambda})}
{\sh^3((b-a)\sqrt{2\lambda})}\\
&&\hskip1.4cm
=\P(I_T\in da, S_T\in db, H_I<H_S)/dadb
\end{eqnarray*}
(cf. Corollary \ref{hit1}). 
Introducing in (\ref{en43}) new variables
via $x=b-a$ and $y=b$ allows us to write 
\begin{equation}
\label{en44} 
\Delta= \int_0^\infty dx\,\varphi(x)\,\left(\Delta_1(x)\Delta_2(x)+\Delta_3(x)\right)
\end{equation} 
where
$$
\Delta_1(x):=\frac{2\lambda\, (\ch(  \beta\sqrt{2\lambda})-1)}{\sh^3(\beta\sqrt{2\lambda})}\,
\exp\left(-(x-\beta)\sqrt{2\la}\ 
 \cth(\beta\,\sqrt{2\la})\right),
$$

\begin{eqnarray*}
&&\hskip-.5cm
\Delta_2(x):={\bf 1
}_{\{\beta<x\}}\int_{x-\beta}^x\sh\left((y-x+\beta)\sqrt{2\lambda}\right)\, dy
=
\frac {\ch(  \beta\sqrt{2\lambda})-1}{\sqrt{{2\lambda}}}\,{\bf 1
}_{\{\beta<x\}},
\end{eqnarray*}
and

\begin{eqnarray*}
&&
\Delta_3(x):= \frac{2\lambda\,\left(\ch(
  x\sqrt{2\lambda})-1\right)}{\sh^3(x\sqrt{2\lambda})}\,
{\bf 1
}_{\{\beta>x\}}\int_{0}^x\sh(y\sqrt{2\lambda})\, dy
\\
&& \hskip1.4cm 
=
 \frac{\sqrt{2\lambda}\,\left(\ch(
  x\sqrt{2\lambda})-1\right)^2}{\sh^3(x\sqrt{2\lambda})}\,
{\bf 1
}_{\{\beta>x\}}.
\end{eqnarray*}
The claimed formula (\ref{en41}) results now easily from (\ref{en44}).
\end{proof}

Before giving results for the joint distribution of $D^+_{T}$ and $D^-_{T}$ we consider the marginal distributions under the condition $ H_I<H_S.$ 

\begin{corollary}
\label{cor41}
For $\al, \beta>0$ 
\begin{equation}
\label{en45}
\P(D^+_{T}\in d\al, H_I<H_S)=
 \frac{\sqrt{2\lambda}\,\left(\ch(
  \al\sqrt{2\lambda})-1\right)^2}{\sh^3(\al\sqrt{2\lambda})}\, d\al,
\end{equation}
\begin{equation}
\label{en46}
\P(D^+_{T}<\al, H_I<H_S)=\frac {\ch(
  \al\sqrt{2\lambda})-1}{2(\ch(\al\sqrt{2\lambda})+1)},
\end{equation}
\begin{equation}
\label{en461}
\P(D^-_{T}<\beta, H_I<H_S)=\frac {(\ch(
  \beta\sqrt{2\lambda})-1)(\ch(
  \beta\sqrt{2\lambda})+2)}{2\,
(\ch(\beta\sqrt{2\lambda})+1)\,\ch(
  \beta\sqrt{2\lambda})}, 
\end{equation}
and
\begin{equation}
\label{en47}
\P(D^+_{T}<\al)=1-\frac 1{\ch(\al\sqrt{2\lambda})}.
\end{equation}
\end{corollary}

\begin{proof}
Letting $\beta\to +\infty$ in (\ref{en41}) yields (\ref{en45}). To
  compute the distribution function in (\ref{en46}) notice that 
\begin{eqnarray*}
&&
\frac{\left(\ch(
  \al\sqrt{2\lambda})-1\right)^2}{\sh^3(\al\sqrt{2\lambda})}
=\frac{\sh(  \al\sqrt{2\lambda}/2)}{2\,\ch^3(
  \al\sqrt{2\lambda}/2)}
\end{eqnarray*}
and hence
\begin{eqnarray*}
&&
\P(D^+_{T}<\al, H_I<H_S)=\frac 12\,\left(1-\frac 1{\ch^2(\al\sqrt{2\lambda}/2)}\right)
\\
&&\hskip4.3cm
=\frac 12\,\left(1-\frac 2{1+\ch(\al\sqrt{2\lambda})}\right),
\end{eqnarray*}
and (\ref{en46}) follows easily. Choosing $\varphi\equiv 1$
in (\ref{en41}), integrating therein and using the above computation 
it is straightforward to derive (\ref{en461}). For (\ref{en47}) use
(\ref{e411}), (\ref{en45}), and (\ref{en461}).
\end{proof}

The proof of Theorem \ref{thm41} shows that the way to express the law of 
$(D^+_{T},D^-_{T})$ in (\ref{en41}) is very natural. However, we need also "more standard" representations obtained when choosing in (\ref{en41}) $\varphi(x)={\bf 1}_{\{ x<\al\}}.$ 

\begin{proposition}
\label{cor42}
The  distribution function of $(D^+_{T},D^-_{T},{\bf 1}_{\{H_I<H_S\}})$ is given
for $\al,\beta\geq 0$ by
\begin{eqnarray}
\label{en48}
&&\hskip-.5cm
\P(D^+_{T}<\al, D^-_{T}<\beta, H_I<H_S)=\frac
12\,\frac{\ch((\al\wedge \beta)\sqrt{2\lambda})-1}{\ch((\al\wedge \beta)\sqrt{2\lambda})+1}
\\
&&\hskip1.5cm\nonumber
+ \frac{\ch(\beta\sqrt{2\lambda})-1}{\ch(\beta\sqrt{2\lambda})(\ch(\beta\sqrt{2\lambda})+1)}
\\
&&\hskip2.5cm\nonumber
\times\,\left(1-\exp\left(-(\al-\beta)\sqrt{2\la}\ 
 \cth(\beta\,\sqrt{2\la})\right)\right)\,{\bf 1}_{\{\beta<\al\}}.
\end{eqnarray}
\end{proposition}
\noindent
Notice that 
\begin{equation}
\label{en485}
H_I<H_S\quad\Leftrightarrow \quad
D^+_{T}\geq D^-_{T},
\end{equation}
explaining the appearance of the indicator function ${\bf 1}_{\{\beta<\al\}}$ in 
(\ref{en48}). Moreover, combining (\ref{en48}) with (\ref{e41}), (\ref{e4105}),
  (\ref{e411}), and  (\ref{en47}) yields

\begin{proposition}
\label{cor43}
The  distribution function of $(D^+_{T},D^-_{T})$ is given
for $\al,\beta> 0$ by
\begin{equation}
\label{en49}
\P(D^+_{T}<\al, D^-_{T}<\beta)=
\begin{cases}
u(\al,\beta),& \al\leq\beta,\\ 
u(\beta,\al),& \al\geq\beta,\\ 
\end{cases}
\end{equation}
where 
$$
u(\al,\beta):=1- 
\frac 1{\ch(\al\sqrt{2\lambda})}- v(\al,\beta)
$$
and 
\begin{equation}
\label{en491}
 v(\al,\beta):=
\frac{(\ch(\al\sqrt{2\lambda})-1)\exp\left(-(\beta-\al)\sqrt{2\la} 
 \cth(\al\,\sqrt{2\la})\right)}
{\ch(\al\sqrt{2\lambda})(\ch(\al\sqrt{2\lambda})+1)}.
\end{equation}
In particular, 
$$
\P(D^+_{T}>\al, D^-_{T}<\beta)=v(\al\wedge\beta,\al\vee\beta).
$$
\end{proposition}
\noindent
We are now able to determine the density function of $(D^+_{T},D^-_{T}).$ Differentiating in (\ref{en49}) leads after some calculations to 
\begin{proposition}
\label{cor44}
For $\al,\beta>0$
\begin{eqnarray}
\label{eq410}
&&
\hskip-1.5cm
\P(D^+_{T}\in d\al\,,\,D^-_{T}\in d\beta)=f_{+,-}(\al\vee\beta,\al\wedge\beta)\,d\al\, d\beta,
\end{eqnarray}
where with $x>y>0$
\begin{eqnarray*}
&&\hskip-1cm
f_{+,-}(x,y) :=
\frac{2\lambda}{\ (\ch(y\sqrt{2\la})+1)^2}
\left(2+\frac{(x-y)\sqrt{2\la}}{\sh(y\,\sqrt{2\la})}\right)
\\
&&
\nonumber
\hskip3cm
\times 
\ \exp\left(-(x-y)\sqrt{2\la}\ \cth(y\,\sqrt{2\la})\right).
\end{eqnarray*}
\end{proposition}

\begin{remark}
To increase understanding of the  distribution of
$(D^+_{T},D^-_{T})$ notice that 
$$
\P(D^+_{T}\in d\al\,,\,D^-_{T}\in d\beta\,|\, H_I<H_S)
=2\,
f_{+,-}(\al,\beta)\,{\bf 1}_{\{\al>\beta\}}\,d\al\, d\beta.
$$
Consequently, for $\al>0,\beta>0$ 
\begin{eqnarray*}
&&\hskip-1cm
\P(D^+_{T}-D^-_{T}\in d\al\,,\,D^-_{T}\in d\beta\,|\, H_I<H_S)
\\
&&\hskip1cm
=
\frac{4\lambda}{\ (\ch(\beta\sqrt{2\la})+1)^2}
\left(2+\frac{\al\sqrt{2\la}}{\sh(\beta\,\sqrt{2\la})}\right)
\\
&&
\hskip2.5cm
\times 
\ \exp\left(-\al\sqrt{2\la}\ \cth(\beta\,\sqrt{2\la})\right)\,d\al\, d\beta.
\end{eqnarray*}
Define two new random variables 
$$
X:=\sqrt{2\la}\left(D^+_{T}-D^-_{T}\right)/\sh(\sqrt{2\la}\,D^-_{T})
$$
and
$$
Y:=\ch(\sqrt{2\la}\,D^-_{T}).
$$
Then for $x>0, y>1$
$$
\P(X\in dx,Y\in dy\,|\, H_I<H_S)= \frac{2(2+x)}{(1+y)^2}\,{\rm e}^{-x\,y}\,
dx\,dy.
$$
In particular, for $y>1$
$$
\P(Y\in dy\,|\, H_I<H_S)= \frac{2(2y+1)}{(y(1+y))^2}\,dy
$$
and for $u>1$
$$
\P(Y>u\,|\, H_I<H_S)= \frac{2}{u(1+u)}.
$$
\end{remark}

Finally, as an application of Proposition \ref{cor44}, we compute the
covariance of $D^+_t$ and $D^-_t$ (up to a fixed time $t$) and determine therefrom 
their correlation coefficient. 
\begin{corollary} 
\label{cor45}
For all $t\geq 0$
$$
\E(D^+_t)=\sqrt{\frac 2\pi}\,\beta(1)\,\sqrt{ t}=\sqrt{\frac \pi 2}
\,\sqrt{ t}\simeq 1.25331\,\sqrt{ t},
$$
$$
\E\left((D^+_t)^2\right)
=
2\,\beta(2)\,t
\simeq 1.83193\, t,
$$
and
$$
\E(D^+_t\,D^-_t)=
\left(2\,\beta(2)-2\log2+1\right)\, t\simeq 1.44564\, t,
$$
where
$$
\beta(n):= \sum_{k=0}^\infty(-1)^k\,(2k+1)^{-n},\quad n=1,2,\dots.
$$
\end{corollary}
\begin{remark}
\label{re00}{\bf 1.} {\sl Recall from Abramowitz and Stegun \cite{abramowitzstegun70}
p. 807, that $\beta(n)$ is called Dirichlet's $\beta$-function. In particular,
$\beta(1)=\pi/2$ and  $\beta(2)=0.91596...$ is Catalan's constant. 
}\hfill\break\hfill
{\bf 2.} {\sl The variance of $D^+$ is 
$$
{\bf Var}(D^+_t):=
\E\left((D^+_t)^2\right)-\left(\E(D^+_t)\right)^2\simeq 0.26113\,t.
$$
The correlation coefficient $\rho$ between $D^+_t$ and $D^-_t$ does not
depend on $t$  and is given by 
\begin{equation}
\label{rho}
\rho:=\frac{\E(D^+_t\,D^-_t)-(\E(D^+_t))^2}{{\bf Var}(D^+_t)}
\simeq 
-0.47936.
\end{equation}
}
\end{remark}

\begin{proof}
From the scaling property of BM it follows that 
\begin{equation}
\label{ddt}
(D^+_T,D^-_T)\ {\mathop=^{\rm{(d)}}}\ \sqrt{T}\,(D^+_1,D^-_1).
\end{equation} 
Since $\E(\sqrt{T})=\sqrt{\pi/\lambda}$ 
and $T$ is independent of the underlying BM we have using (\ref{en47})
\begin{eqnarray*}
&&\hskip-1.5cm
\E(D^+_1)=\frac{\E(D^+_T)}{\E(\sqrt T)}=\sqrt{\frac\lambda\pi}\,\int_0^\infty 
\P(D^+_{T}>\al)\, d\al
\\
&&\hskip1.8cm
=(\sqrt{2\pi})^{-1}\,\int_0^\infty\frac 1{\ch u}\, du=
\sqrt{\frac 2\pi}\, \beta(1),
\end{eqnarray*}
where the series expansion 
$$
(\ch u)^{-1}=2\,{\rm e}^{-u}\,\sum_{k=0}^\infty (-1)^k\, {\rm e}^{-2uk}
$$ 
is used. This yields the first formula in the statement of the corollary. For the second formula we compute similarly:
\begin{eqnarray*}
&&
\E((D^+_1)^2)
=\frac{\E((D^+_T)^2)}{\E(T)}
=
2\,\lambda\,\int_0^\infty 
\P(D^+_{T}>\al)\,\al\, d\al
\\
&&
\hskip4.2cm
=\int_0^\infty\frac u{\ch u}\, du
=
 2\,\beta(2).
\end{eqnarray*}
Next we determine $\E(D^+_t\,D^-_t).$ Firstly,
\begin{eqnarray*}
&&\hskip-1.5cm
\E(D^+_1\,D^-_1)=\frac{\E(D^+_T\,D^-_T)}{\E(T)} 
 =\lambda\,\int_0^\infty\int_0^\infty 
\P(D^+_{T}>\al,D^-_{T}>\beta)\,d\al\,d\beta.
\end{eqnarray*}
Supposing $\al<\beta$ we have from (\ref{en49})
\begin{eqnarray*}
&&\hskip-1.5cm
\P(D^+_{T}<\al,D^-_{T}<\beta)
=1- 
\frac 1{\ch(\al\sqrt{2\lambda})}- v(\al,\beta)\\
&&\hskip2.3cm
=\P(D^+_{T}<\al)-v(\al,\beta)
\end{eqnarray*}
with $v$ as given in (\ref{en491}).
Consequently,
$$
\P(D^+_{T}<\al,D^-_{T}>\beta)= v(\al,\beta),
$$
and
$$
\P(D^-_{T}>\beta) -\P(D^+_{T}>\al,D^-_{T}>\beta)= v(\al,\beta),
$$
which gives 
$$
\P(D^+_{T}>\al,D^-_{T}>\beta)= 
\frac 1{\ch(\beta\sqrt{2\lambda})}- v(\al,\beta),\qquad \al<\beta.
$$
Hence it holds
$$
\int_\al ^\infty \P(D^+_{T}>\al,D^-_{T}>\beta)\, d\beta= 
\int_\al ^\infty \frac 1{\ch(\beta\sqrt{2\lambda})}\, d\beta
- \int_\al ^\infty v(\al,\beta)\, d\beta,
$$
where, by elementary integration, 
$$
\int_\al ^\infty v(\al,\beta)\, d\beta=
\frac{(\ch(\al\sqrt{2\lambda})-1)\, \sh(\al \sqrt{2\la})} 
{\sqrt{2\la}\,(\ch(\al\sqrt{2\lambda}))^2(\ch(\al\sqrt{2\lambda})+1)}.
$$
By symmetry, we have 
\begin{eqnarray*}
&&\hskip-.5cm
\E(D^+_1\,D^-_1)
=2\lambda\left(\int_0^\infty d\al \int_\al ^\infty d\beta \frac 1{\ch(\beta\sqrt{2\lambda})}
- \int_0^\infty d\al \int_\al ^\infty d\beta\ v(\al,\beta)\right)
\\
&&\hskip1.5cm
=\int_0^\infty \frac{u}{\ch u}\,du-\int_0^\infty \frac{(\ch u-1)\sh u}{\ch^2 u(\ch u+1)}\,du
\\
&&\hskip1.5cm
=2\sum_{k=0}^\infty (-1)^k(2k+1)^{-2}-\int_1^\infty \frac{x-1}{x^2(x+1)}\, dx
\\
&&\hskip1.5cm
=2\beta(2)  -2\log 2-1.
\end{eqnarray*}
\end{proof}
\begin{remark}
Recall that
$$
\E( T^{\,p/2})=\lambda\int_0^\infty t^{p/2}\,{\rm e}^{-\lambda t}\, dt= 
{\Gamma((2+p)/2)}/{\lambda^{p/2}}.
$$
With this formula it is now fairly straightforward to connect the $p$th moment of $D^+_1$ with Dirichlet's $\beta$-function. Indeed, for $p\geq 1$ we have
\begin{eqnarray*}
&&\hskip-.5cm
\E((D^+_1)^p)=\frac{\E((D^+_T)^p)}{\E( T^{p/2})}=\frac{\lambda^{p/2}\,p}{\Gamma(\frac{2+p}{2})}\ \int_0^\infty 
\P(D^+_{T}>\al)\, \al^{p-1}\,d\al
\\
&&\hskip3.6cm
=\frac{p}{2^{p/2}\,\Gamma(\frac{2+p}{2})}\ \int_0^\infty  \frac{t^{p-1}\, dt}{\ch t}
\\
&&\hskip3.6cm
=\frac{2\,p\,\Gamma(p)}{2^{p/2}\,\Gamma(\frac{2+p}{2})}\ \beta(p).
\end{eqnarray*}
\end{remark}

\noindent
{\bf Acknowledgement.} We thank Gabor Szekely for posing the problem
of finding the covariance of maximum increase and decrease of BM,
and Margret Halldorsdottir for observing the connection with
Catalan's constant.

\end{document}